# You Can Fool Some People Sometimes...

Rasa Karapandza♣ and Milos Bozovic♠

Universitat Pompeu Fabra and SECCF

**Abstract**

We develop an empirical procedure to qunatify future company performance based on top management promises. We find that the number of future tense sentence occurrences in 10-K reports is significantly negatively correlated with the return as well as with the excess return on the company stock price. We extrapolate the same methodology to US presidential campaigns since 1960 and come to some startling conclusions.

## 1. Introduction

Annual reports, including 10-K reports in the United States, have two major functions: a) to inform current and potential investors, as well as the regulators, about significant developments that shaped previous year's company results and b) to help investors forecast company results for the next period. In such a situation it is conceivable that the top management of companies may use these reports to make excessive promises in an attempt to convince the investors that the situation is *rosy* even when it is not. To get an idea about a potential relation between such promises and the subsequent company performance we looked at 555 10-K reports of the American S&P 100 companies during the period from 1993 to 2003. We found a whopping 313,852 times that future tense constructions were used in the reports (we define these as sentences that contain either *will*, *shall*, or *going to*). Importantly, companies who used less future tense sentences in their reports achieved significantly higher returns and excess returns in the subsequent fiscal year.

A similar pattern of behavior can be observed among politicians and can give us a simple way of predicting the outcome of the American presidential elections! Namely, in all presidential debates in the period 1960–2000 for which the transcripts from the debates were available (in 1984 two debate transcripts were not available, while the debates were not conducted in 1964, 1968 and 1972), a presidential candidate who made less use of future tense sentences won the popular vote whether he was an incumbent or a

We would like to thank our professors Jaume Garcia, Christian Hafke, Jose Marin, Branko Urosevic and our coleagues Gueorgui Kolev and Blaz Zakelj for their guidance and unselfish help. All remaining errors in this paper are our own. Both authors are PhD students at Department of Economics, Universitat Pompeu Fabra, Ramon Trias Fargas 25–27, 08005 Barcelona, Spain (*www.upf.edu*) and are affiliated with South European Center for Contemporary Finance, Faculty of Economics, Kamenicka 6, 11000 Belgrade, Serbia-Montenegro (*www.seccf.org*).

♣ Corresponding author: Department of Economics, Universitat Pompeu Fabra, Ramon Trias Fargas 25-27, 08005 Barcelona, Spain Tel.: (34) 679 965 103; Email: *rasa.karapandza@upf.edu*

♠ Email: *milos.bozovic@upf.edu*

challenger. According to this simple rule of thumb, this year's winner of the popular vote should be President George W. Bush.

The rest of the paper is organized as follows. Section 2 gives a short overview of the literature. Section 3 describes the methods and the results obtained in the case of companies. Section 4 applies the same basic idea to discuss the signals from presidential debates. Conclusions and future research directions are presented in Section 5. All tables are listed in Section 6.

## 2. Literature overview

Corporate annual reports are an important device for communication between management and shareholders (see Lee and Tweedie, 1975a, 1975b, 1976; Bartlet and Chandler, 1997). It is widely accepted that a major role for financial reporting is to help participants in the financial markets forecast company performance.

Fishman and Hagerty (1989) provide a rationale for voluntary disclosures of information by firms. They model a relationship between investment decisions by firms and the efficiency of the market prices on their securities and show that more efficient security prices can lead to more efficient investment decisions. This provides firms with an incentive to increase price efficiency by voluntarily disclosing information about the firm. Naser and Nuseibeh (2003) study the quality of information in annual reports of Saudi nonfinancial companies and find that Saudi companies disclose information more than the minimum required by law but at the level that is relatively low.

The creation of the World Wide Web provides investors with new opportunities for mining of information about companies. The amount of company information available is substantial. This information is of great interest to company customers, stockholders, creditors, auditors, financial analysts, and management. Ettredge, Richardson, and Scholz (2001) investigate corporate Web site financial disclosure practices. Dull, Graham, and Baldwin (2003) also study the effect of Internet-based financial statements on investment decisions.

Extracting valuable signals from a giant pool of information is a daunting task. Different data mining tools, such as neural networks, have been proposed to help with the task (see Barbro *et al.* (2001), for example). In contrast with the complex data techniques available in the literature, our method of signal extraction from financial statements is quite simple, yet gives strikingly robust predictions about future performance of companies.

## 3. Results

We analyze the disclosure content of 10-K reports for 88 companies included in the S&P 100 portfolio in July 2004. Reports are obtained from the web site of Security and Exchange Commission (SEC) for the period from 1993 to 2003. In total we collected 555 10–K reports, since not all of them were available for all of the companies in the sample (because of the lack of data, only 88 out of 100 companies were used). Data on dividend-adjusted share prices of companies were obtained from Yahoo! Finance.

For a proxy of the informational disclosure content we use the amount of future tense occurrences in the 10-K reports. The variables that we use are:
- The number of occurrences of the verb *will* in 10–K at year $t$ (denoted by $w_t$)
- Return on the company stock in year $t+1$ (denoted by $R_{t+1}$)[1]
- Risk free rate in year $t+1$ (denoted by $r_{ft\ t+1}$)
- The time dummies, *i.e.* the variables that take value 1 for a given year and 0 for all other years.

Throughout the paper, we use panel-data type regressions with random effects. To correct for possible heteroscedasticity and autocorrelation within companies we use pooled data regressions utilizing the robust variance estimator (see, *e.g.*, Wooldrige, 2002).

We first regress the number of *wills* on the number of *shalls* and find significant positive correlation between the two, with $R^2=0.56$. (Sensibly, a well-crafted report has relatively balanced amount of *wills* and *shalls*). Furthermore, in all 555 reports we find only 53 occurrences of the phrase *going to* (in comparison with 313,852 of *wills* and *shalls*) and thus discard that phrase from consideration (we do not discard it in case of the presidential debates, see below). Therefore, as a proxy for the total number of future tense sentences we use only $w_t$.

In order to investigate the determinants of the return on a company stock in year $t+1$ we commence with a panel-data random-effects type regression. The dependent variable is the stock return in year $t+1$; the explanatory variables are: the number of *will* occurrences at year $t$ reports, controlled for economy-wide shocks by including a full set of time dummies (the omitted category is year 1993). In other words, we run

$$R_{t+1} = a + bw_t + d_{1994}\ dummy(1994) + \ldots + d_{2003}\ dummy(2003) + u_t$$

with random effects. The results are summarized in Table 1. Notice that the estimator for $b$ is negative and significant at the 95% confidence level.

Next, to control for heteroscedasticity and any form of autocorrelation within companies we, next, run a regression on the pooled data using the robust variance estimators. Results are given in Table 2. It can be seen that the coefficient associated with $b$ is still negative and significant at the 99% confidence level.

---

[1] Actually, we use $R_{t+1} = \ln S_{t+1} - \ln S_t$, where $S_t$ is the stock price at time $t$, consistently throughout the paper.

Similar regressions are performed for excess stock return. The results of a panel data random-effects regression are presented in Table 3. Again, we control for economy wide effects and include full set of time dummy variables. The specification follows:

$$R_{t+1} - r_{f\,t+1} = a + bw_t + d_{1994}\,dummy(1994) + \ldots + d_{2003}\,dummy(2003) + u_t$$

As before, $b$ is negative and significant at 95% confidence level. The same specification estimated on the pooled data using the robust variance estimator yields negative and, again, highly significant coefficient (99% significance level). Results are presented in Table 4.

Note that the results of the last two regressions could be expected since $r_{f\,t}$ is economy-wide variable and varies only across time (that is, it is not company-specific).

## 4. U.S. Presidential Debates

Voters in the U.S. presidential elections seem to prefer candidates who offer fewer promises during the pre-election debates! We demonstrate this by analyzing the occurrence of future tense sentences (as a proxy for promises) in all of the public presidential debates from 1960 to 2004, whenever transcripts of the debates were available. We obtained transcripts from the web site of the Commission on Presidential Debates (*www.debates.org*). The results are given in Figure 1. Every time the popular vote was won by the candidate, whether he was an incumbent or a challenger, that used the smaller overall number of either *will*, *shall*, or *going to* constructions. The *t*-statistics of the difference between the counts for Loser and Winner is significant at 90% confidence level. In particular, the algorithm shows that Al Gore won the popular votes in the 2000 presidential elections (he lost subsequently due to the Electoral College rule).

**Figure 1**

Future tense usage comparison in U.S. presidential debates. By *winner* we denote the candidate that won majority of popular votes.

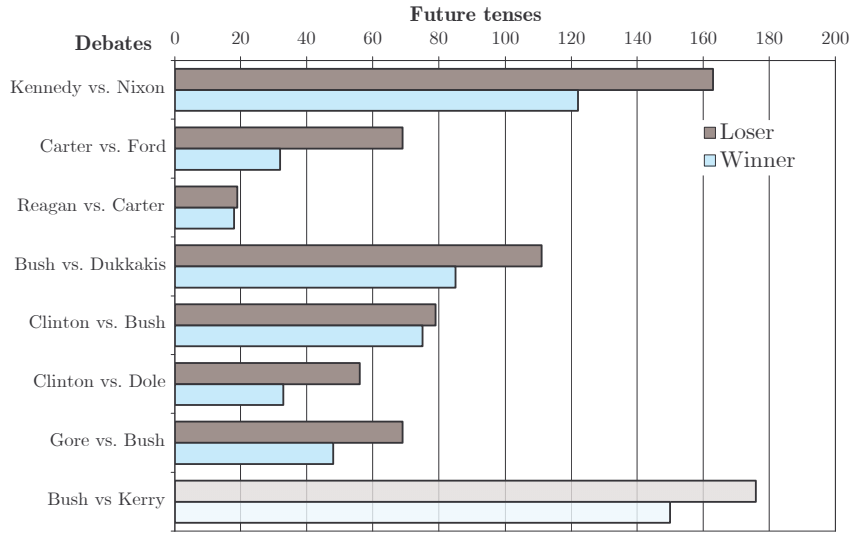

In the 2004 presidential debates John F. Kerry used 176 *wills, shalls* or *going-tos,* while George W. Bush used them only 150 times. Based on this alone, one would be drawn to conclude that George W. Bush will win majority of popular votes in the 2004 presidential elections.

Of course, we are aware of many problems that such an analysis may encounter. The key among them is a relatively small number of televised presidential debates thus far. However, we believe that the extension of our analysis to debates prior to 1960 would not make much sense, since their influence on the voting population did not have the same direct impact - they were not televised and thus were reaching a smaller fraction of the population.

## 5. Conclusion

In this paper we have demonstrated how a simple empirical procedure can provide a method for signal extraction out of the 10-K reports. According to our preliminary results there exists a highly significant negative correlation between company returns and the number of future tense sentences used in their 10-K reports. Similarly, there is a significant negative correlation between the excess return on a company's stock and the number of future tense sentences they use in the 10-K reports.

A similar negative correlation between the number of occurrences of future tense sentences in the U.S. presidential debates on the results of the popular votes in the corresponding presidential election is observed. In particular, thus far, the candidate that won the popular votes in an election has always been the one who used fewer future-tense constructions. Based on this simple rule of thumb alone, the algorithm predicts that the winner of the popular votes in the 2004 presidential elections will be George W. Bush.

Let us conclude by reminding you of the famous statement by Abraham Lincoln (1809–1865):

*You can fool all the people some of the time, and some of the people all the time, but you cannot fool all the people all the time.*

# 6. Tables

## Table 1

Results of a panel data random effects type of regression of the return of a company in year $t+1$ on a number of occurrences of the verb *will* in 10–K reports in year $t$ (denoted by $w_t$), and full set of dummies excluding the first year in the dataset (1993), since we include constant term in the regression:

| Random-effects GLS regression | | | | Number of obs | 555 |
|---|---|---|---|---|---|
| Group variable (i): id | | | | Number of groups | 88 |
| R-sq: | within | 0.2400 | Obs per group: | min | 1 |
| | between | 0.1501 | | avg | 6.3 |
| | overall | 0.2212 | | max | 10 |
| Random effects u_i ≈ Gaussian | | | | Wald chi2(1) | 154.27 |
| corr(u_i, X) = 0 (assumed) | | | | Prob > chi2 | 0.0000 |

| Dep var. $R_{t+1}$ | Coef. | Std. Err. | z | P>|z| | [95% Conf. Interval] | |
|---|---|---|---|---|---|---|
| $w_t$ | −.0001965 | .0000833 | −2.36 | 0.018 | −.0003598 | −.0000332 |
| dum1994 | −.0494584 | .2474785 | −0.20 | 0.842 | −.5345073 | .4355905 |
| dum1995 | .2101371 | .2477426 | 0.85 | 0.396 | −.2754294 | .6957036 |
| dum1996 | .1037145 | .2477875 | 0.42 | 0.676 | −.3819401 | .589369 |
| dum1997 | .1309552 | .2476234 | 0.53 | 0.597 | −.3543777 | .6162881 |
| dum1998 | .162265 | .2476778 | 0.66 | 0.512 | −.3231747 | .6477046 |
| dum1999 | .1397685 | .2474706 | 0.56 | 0.572 | −.3452649 | .624802 |
| dum2000 | .0699918 | .2474654 | 0.28 | 0.777 | −.4150314 | .555015 |
| dum2001 | −.164591 | .2475036 | −0.67 | 0.506 | −.6496893 | .3205072 |
| dum2002 | −.1669186 | .2473577 | −0.67 | 0.500 | −.6517308 | .3178935 |
| dum2003 | .0533142 | .2467967 | 0.22 | 0.829 | −.4303984 | .5370267 |
| cons | .1230746 | .2450847 | 0.50 | 0.616 | −.3572825 | .6034317 |

**Table 2**

Results of a regression on the pooled data using the robust variance estimator allowing for heteroscedasticity and autocorrelation within companies of return of a company in year $t+1$ and number of occurrences of the verb *will* in 10–K reports in year $t$ (denoted by $w_t$), and dummy variables for every year excluding the first year in the dataset (1993), since we include constant term in the regression:

| Regression with robust standard errors | | | | | Number of obs | 555 |
|---|---|---|---|---|---|---|
| Number of clusters (id) = 88 | | | | | R-squared | 0.2212 |
| | | Robust | | | | |
| Dep var. $R_{t+1}$ | Coef. | Std. Err. | z | P>\|z\| | [95% Conf. Interval] | |
| $w_t$ | −.0001965 | .0000734 | −2.68 | 0.009 | −.0003423 | −.0000506 |
| dum1994 | −.0494584 | .0256026 | −1.93 | 0.057 | −.1003463 | .0014295 |
| dum1995 | .2101371 | .0304549 | 6.90 | 0.000 | .1496048 | .2706694 |
| dum1996 | .1037145 | .0316055 | 3.28 | 0.001 | .0408952 | .1665338 |
| dum1997 | .1309552 | .0374868 | 3.49 | 0.001 | .056446 | .2054643 |
| dum1998 | .162265 | .0428301 | 3.79 | 0.000 | .0771355 | .2473944 |
| dum1999 | .1397685 | .0393589 | 3.55 | 0.001 | .0615384 | .2179986 |
| dum2000 | .0699918 | .0469311 | 1.49 | 0.139 | −.0232889 | .1632725 |
| dum2001 | −.164591 | .0349906 | −4.70 | 0.000 | −.2341387 | −.0950434 |
| dum2002 | −.1669186 | .0280812 | −5.94 | 0.000 | −.2227331 | −.1111042 |
| dum2003 | .0533142 | .0207075 | 2.57 | 0.012 | .0121558 | .0944725 |
| cons | .1230746 | .0013208 | 93.18 | 0.000 | .1204494 | .1256998 |

## Table 3

Results of a panel data random effects type of regression of the excess return of a company in year $t+1$ on a number of occurrences of the verb *will* in 10–K reports in year $t$ (denoted by $w_t$), and dummy variables for every year excluding the first year in the dataset (1993), since we include constant term in the regression:

| Random-effects GLS regression | | | | Number of obs | | 555 |
|---|---|---|---|---|---|---|
| Group variable (i): id | | | | Number of groups | | 88 |
| R-sq: | within | 0.2171 | Obs per group: | | min | 1 |
| | between | 0.1303 | | | avg | 6.3 |
| | overall | 0.1989 | | | max | 10 |
| Random effects u_i ≈ Gaussian | | | | Wald chi2(1) | | 134.79 |
| corr(u_i, X) = 0 (assumed) | | | | Prob > chi2 | | 0.0000 |
| Dep. var. $R_{t+1} - rf_{t+1}$ | Coef. | Std. Err. | z | P>\|z\| | [95% Conf. Interval] | |
| $w_t$ | −.0001965 | .0000833 | −2.36 | 0.018 | −.0003598 | −.0000332 |
| dum1994 | −.0498584 | .2474785 | −0.20 | 0.840 | −.5349073 | .4351904 |
| dum1995 | .1746371 | .2477426 | 0.70 | 0.481 | −.3109294 | .6602035 |
| dum1996 | .0878145 | .2477875 | 0.35 | 0.723 | −.3978401 | .573469 |
| dum1997 | .1098552 | .2476234 | 0.44 | 0.657 | −.3754777 | .5951881 |
| dum1998 | .144865 | .2476778 | 0.58 | 0.559 | −.3405747 | .6303046 |
| dum1999 | .1296685 | .2474706 | 0.52 | 0.600 | −.355365 | .614702 |
| dum2000 | .0437918 | .2474654 | 0.18 | 0.860 | −.4412314 | .528815 |
| dum2001 | −.177691 | .2475036 | −0.72 | 0.473 | −.6627893 | .3074072 |
| dum2002 | −.1535187 | .2473577 | −0.62 | 0.535 | −.6383308 | .3312935 |
| dum2003 | .0747141 | .2467967 | 0.30 | 0.762 | −.4089984 | .5584267 |
| cons | .0880746 | .2450847 | 0.36 | 0.719 | −.3922825 | .5684317 |

## Table 4

Results of a regression on the pooled data using the robust variance estimator allowing for heterocedasticity and autocorrelation within companies of the excess return of a company in year $t+1$ on a number of occurrences of the verb *will* in 10–K reports in year $t$ (denoted by $w_t$), and dummy variables for every year excluding the first year in the dataset (1993), since we include constant term in the regression:

| Regression with robust standard errors | | | | | Number of obs | 555 |
|---|---|---|---|---|---|---|
| Number of clusters (id) = 88 | | | | | R-squared | 0.1989 |
| Dep var. $R_{t+1} - rf_{t+1}$ | Coef. | Robust Std. Err. | z | P>\|z\| | [95% Conf. Interval] | |
| wt | −.0001965 | .0000734 | −2.68 | 0.009 | −.0003423 | −.0000506 |
| dum1994 | −.0498584 | .0256026 | −1.95 | 0.055 | −.1007463 | .0010295 |
| dum1995 | .1746371 | .0304549 | 5.73 | 0.000 | .1141047 | .2351694 |
| dum1996 | .0878145 | .0316055 | 2.78 | 0.007 | .0249952 | .1506338 |
| dum1997 | .1098552 | .0374868 | 2.93 | 0.004 | .035346 | .1843643 |
| dum1998 | .144865 | .0428301 | 3.38 | 0.001 | .0597355 | .2299944 |
| dum1999 | .1296685 | .0393589 | 3.29 | 0.001 | .0514384 | .2078986 |
| dum2000 | .0437918 | .0469311 | 0.93 | 0.353 | −.0494889 | .1370725 |
| dum2001 | −.177691 | .0349906 | −5.08 | 0.000 | −.2472387 | −.1081434 |
| dum2002 | −.1535187 | .0280812 | −5.47 | 0.000 | −.2093331 | −.0977042 |
| dum2003 | .0747141 | .0207075 | 3.61 | 0.001 | .0335558 | .1158725 |
| cons | .0880746 | .0013208 | 66.68 | 0.000 | .0854494 | .0906998 |